\documentclass[12pt,twoside]{article}
\usepackage[english]{babel}
\usepackage[latin1]{inputenc}
\usepackage{amsmath}
\usepackage{amssymb,amsfonts}
\usepackage{graphicx}
\usepackage{times,amssymb,amscd}

\newcommand{\bC}{\mathbf{C}}
\newcommand{\bE}{\mathbf{E}}
\newcommand{\bG}{\mathbf{G}}
\newcommand{\bH}{\mathbf{H}}

\newcommand{\bL}{\mathbf{L}}

\newcommand{\bR}{\mathbf{R}}
\newcommand{\bS}{\mathbf{S}}

\newcommand{\bT}{\mathbf{T}}

\newcommand{\bt}{\mathbf{t}}

\newcommand{\cP}{\mathcal{P}}
\newcommand{\cS}{\mathcal{S}}

\newcommand{\cC}{\mathcal{C}}
\newcommand{\cH}{\mathcal{H}}

\newcommand{\EUC}{\bE^3}

\newcommand{\SXR}{\bS^2\!\times\!\bR}
\newcommand{\HXR}{\bH^2\!\times\!\bR}
\newcommand{\SLR}{\widetilde{\bS\bL_2\bR}}
\newcommand{\NIL}{\mathbf{Nil}}
\newcommand{\SOL}{\mathbf{Sol}}

\newtheorem{Definition}{Definition}[section]
\newtheorem{Remark}{Remark}[section]

\newtheorem{Theorem}{Theorem}[section]
\begin{document}
\pagestyle{myheadings}
\markboth{\centerline{Jen\H o Szirmai}}
{Pappus's hexagon theorem for geodesic curves...}
\title
{Pappus's hexagon theorem for geodesic curves in $\SLR$ geometry
}

\author{Jen\H o Szirmai \\
\normalsize Department of Algebra and Geometry, Institute of Mathematics,\\
\normalsize Budapest University of Technology and Economics, \\
\normalsize M\H uegyetem rkp. 3., H-1111 Budapest, Hungary \\
\normalsize szirmai@math.bme.hu \\
\date{\normalsize{\today}}}


\maketitle
\begin{abstract}
In this paper, we continue our investigation of the structure 
of Thurston geometries with non-constant curvature. 
We generalize the notion of cylinders, prove that they contain the types of geodesic curves.
Moreover, we formulate and prove a 
statement analogous to Pappus's hexagon theorem for $\SLR$ geometry. 
\end{abstract}

\newtheorem{theorem}{Theorem}[section]
\newtheorem{corollary}[theorem]{Corollary}
\newtheorem{conjecture}[theorem]{Conjecture}
\newtheorem{lemma}[theorem]{Lemma}
\newenvironment{example}{Example}



\section{Introduction}
Pappus's hexagon theorem plays a fundamental role in the axiomatic construction 
of projective geometry. As seen in \cite{Sz26}, this theorem can be extended to 
$\NIL$ space, where geodesics -- within the unified projective model of Thurston geometries 
(see \cite{M97}) -- are helix-like curves lying on $\NIL$ cylinders 
(which are cylinders in the Euclidean sense as well). 
The $\NIL$ version of the theorem highlights the fact that, unlike in classical geometries, 
the geodesics do not belong to a single isometry class; consequently, 
Pappus's hexagon theorem is modified accordingly.

\begin{Remark}
\begin{enumerate}
\item Pappus's hexagon theorem is well known in the classical geometries of 
constant curvature. In the $\NIL$ geometry 
(see \cite{Sz26}) we extended it for cylinders and helix-like geodesic curves. 
We are working on an extension for the $\SXR$ and $\HXR$ geometries. 

As we demonstrate in this article the Pappus's theorem is extendable 
for geodesics of $\SLR$ geometry. 

Among the Thurston geometries, $\SOL$ is the only one for which 
this cannot be stated, as the equations 
for its geodesic curves are known exactly only in special cases. (On my opinion it will also be extendable to $\SOL$ geometry.)
\item We examined several interesting problems related to the Thurston geometries 
with non-constant curvature in papers \cite{Cs23}, \cite{CsSz16}, \cite{CsSz23}, \cite{MSz12}, \cite{MSzV}, 
\cite{Sz07}, \cite{Sz13-1}, \cite {Sz14}, \cite{Sz13-2}, \cite{Sz23-2}, \cite{Sz24}. 
\end{enumerate}
\end{Remark}

{In this paper, we continue the investigation of similar questions regarding {\it geodesic curves} in $\SLR$ geometry. 
We generalize the notion of $\SLR$ cylinders where the notion of circular cylinders
is introduced in \cite{Sz24}. We prove 
that the geodesic curves lie on generalized $\SLR$ cylinders (see Definition 3.1). Moreover, we prove that the radius of the 
generalized cylinders 
characterizes the equivalence classes of $\SLR$ geodesics. 
The main objective of this work to analysing cylinder-like Pappus's theorem and 
formulate a statement analogous to Pappus's theorem for 
generalized $\SLR$ cylinders and geodesic curves.}
\section{Preliminaries and background}
The real $ 2\times 2$ matrices $\begin{pmatrix}
         d&b \\
         c&a \\
         \end{pmatrix}$ with unit determinant $ad-bc=1$
constitute a Lie transformation group by the usual product operation, taken to act on row matrices as on point coordinates on the right as follows
\begin{equation}
\begin{gathered}
(z^0,z^1)\begin{pmatrix}
         d&b \\
         c&a \\
         \end{pmatrix}=(z^0d+z^1c, z^0 b+z^1a)=(w^0,w^1)\\
\mathrm{with} \ w=\frac{w^1}{w^0}=\frac{b+\frac{z^1}{z^0}a}{d+\frac{z^1}{z^0}c}=\frac{b+za}{d+zc}, \tag{2.1}
\end{gathered}
\end{equation}
as action on the complex projective line $\bC^{\infty}$ (see \cite{M97}, \cite{MSz}).
This group is a $3$-dimensional manifold, because of its $3$ independent real coordinates and with its usual neighbourhood topology (\cite{R}, \cite{S}, \cite{T}).
In order to model the above structure in the projective sphere $\cP \cS^3$ and in the projective space $\cP^3$ (see \cite{M97}),
we introduce the new projective coordinates $(x^0;x^1;x^2;x^3)$ where
$a:=x^0+x^3, \ b:=x^1+x^2, \ c:=-x^1+x^2, \ d:=x^0-x^3$
with the positive, then the non-zero multiplicative equivalence as projective freedom in $\cP \cS^3$ and in $\cP^3$, respectively.
Then it follows that $0>bc-ad=-x^0x^0-x^1x^1+x^2x^2+x^3x^3$
describes the interior of the above one-sheeted hyperboloid solid $\cH$ in the usual Euclidean coordinate simplex with the origin
$E_0(1;0;0;0)$ and the ideal points of the axes $E_1^\infty(0;1;0;0)$, $E_2^\infty(0;0;1;0)$, $E_3^\infty(0;0;0;1)$.
We consider the collineation group ${\bf G}_*$ that acts on the projective sphere $\cS\cP^3$  and preserves a polarity i.e. a scalar product of signature
$(- - + +)$, this group leaves the one sheeted hyperboloid solid $\cH$ invariant.
We have to choose an appropriate subgroup $\mathbf{G}$ of $\mathbf{G}_*$ as isometry group, then the universal covering group and space
$\widetilde{\cH}$ of $\cH$ will be the hyperboloid model of $\SLR$ \cite{M97}.

The specific isometries $\bS(\phi)$ $(\phi \in \bR )$ constitute a one parameter group given by the matrices:
\begin{equation}
\begin{gathered} \bS(\phi):~(s_i^j(\phi))=
\begin{pmatrix}
\cos{\phi}&\sin{\phi}&0&0 \\
-\sin{\phi}&\cos{\phi}&0&0 \\
0&0&\cos{\phi}&-\sin{\phi} \\
0&0&\sin{\phi}&\cos{\phi}
\end{pmatrix}
\end{gathered} \tag{2.2}
\end{equation}
The elements of $\bS(\phi)$ are the so-called {\it fibre translations}. We obtain a unique fibre line to each $X=(x^0;x^1;x^2;x^3) \in \widetilde{\cH}$
as the orbit under the right action of $\bS(\phi)$ on $X$. The coordinates of points lying on the fibre line through $X$ can be expressed
as the images of $X$ by $\bS(\phi)$:
\begin{equation}
\begin{gathered}
(x^0;x^1;x^2;x^3) \stackrel{\bS(\phi)}{\longrightarrow} {(x^0 \cos{\phi}-x^1 \sin{\phi}; x^0 \sin{\phi} + x^1 \cos{\phi};} \\ {x^2 \cos{\phi} + x^3 \sin{\phi};-x^2 \sin{\phi}+
x^3 \cos{\phi})}.
\end{gathered} \tag{2.3}
\end{equation}
The points of a fibre line through
$X$ by usual inhomogeneous Euclidean coordinates $x=\frac{x^1}{x^0}$, $y=\frac{x^2}{x^0}$, $z=\frac{x^3}{x^0}$, $x^0\ne 0$ are given by
\begin{equation}
\begin{gathered}
(1;x;y;z) \stackrel{\bS(\phi)}{\longrightarrow} {\Big( 1; \frac{x+\tan{\phi}}{1-x \tan{\phi}}; \frac{y+z \tan{\phi}}{1-x \tan{\phi}};
\frac{z - y \tan{\phi}}{1-x \tan{\phi}}\Big)}
\end{gathered} \tag{2.4}
\end{equation}
for the projective space $\cP^3$, where ideal points (at infinity) conventionally occur.

In (2.3) and (2.4) we can see the $2\pi$ periodicity of $\phi$, moreover the (logical) extension to $\phi \in \bR$, as real parameter, to have
the universal covers $\widetilde{\cH}$ and $\SLR$, respectively, through the projective sphere $\cP\cS^3$. The elements of the isometry group of
$\mathbf{SL_2R}$ (and so by the above extension the isometries of $\SLR$) can be described by the matrix $(a_i^j)$ (see \cite{M97} and \cite{MSz})
Moreover, we have the projective proportionality, of course.
We define the {\it translation group} $\bG_T$, as a subgroup of the isometry group of $\mathbf{SL_2R}$,
the isometries acting transitively on the points of ${\cH}$ and by the above extension on the points of $\SLR$ and $\widetilde{\cH}$.
$\bG_T$ maps the origin $E_0(1;0;0;0)$ onto $X(x^0;x^1;x^2;x^3)$. These isometries and their inverses (up to a positive determinant factor)
are given by the following matrices:
\begin{equation}
\begin{gathered} \bT:~(t_i^j)=
\begin{pmatrix}
x^0&x^1&x^2&x^3 \\
-x^1&x^0&x^3&-x^2 \\
x^2&x^3&x^0&x^1 \\
x^3&-x^2&-x^1&x^0
\end{pmatrix}.
\end{gathered} \tag{2.5}
\end{equation}
The rotation about the fibre line through the origin $E_0(1;0;0;0)$ by angle $\omega$ $(-\pi<\omega\le \pi)$ can be expressed by the following matrix
(see \cite{M97})
\begin{equation}
\begin{gathered} \bR_{E_O}(\omega):~(r_i^j(E_0,\omega))=
\begin{pmatrix}
1&0&0&0 \\
0&1&0&0 \\
0&0&\cos{\omega}&\sin{\omega} \\
0&0&-\sin{\omega}&\cos{\omega}
\end{pmatrix},
\end{gathered} \tag{2.6}
\end{equation}
and the rotation $\bR_X(\omega)$ about the fibre line through $X(x^0;x^1;x^2;x^3)$ by angle $\omega$ can be derived by formulas (2.5) and (2.6):
\begin{equation}
\bR_X(\omega)=\bT^{-1} \bR_{E_O} (\omega) \bT:~(r_i^j(X,\omega)).
\tag{2.7}
\end{equation}
Horizontal intersection of the hyperboloid solid $\cH$ with the plane $E_0 E_2^\infty E_3^\infty$ provides the
{\it hyperbolic $\mathbf{H}^2$ base plane} of the model $\widetilde{\cH}=\SLR$.
The fibre through $X$ intersects the base plane $z^1=x=0$ in the foot point
\begin{equation}
\begin{gathered}
Z=(z^0=x^0 x^0+x^1x^1; z^1=0; z^2=x^0x^2-x^1x^3;z^3=x^0x^3+x^1x^2).
\end{gathered} \tag{2.8}
\end{equation}
After \cite{M97}, we introduce the so-called hyperboloid parametrization as follows
\begin{equation}
\begin{gathered}
x^0=\cosh{r} \cos{\phi}, ~ ~
x^1=\cosh{r} \sin{\phi}, \\
x^2=\sinh{r} \cos{(\theta-\phi)}, ~ ~
x^3=\sinh{r} \sin{(\theta-\phi)},
\end{gathered} \tag{2.9}
\end{equation}
where $(r,\theta)$ are the polar coordinates of the base plane and $\phi$ is just the fibre coordinate. We note that
$$-x^0x^0-x^1x^1+x^2x^2+x^3x^3=-\cosh^2{r}+\sinh^2{r}=-1<0.$$
The inhomogeneous coordinates corresponding to (2.9), that play an important role in the later visualization of prism tilings in $\EUC$,
are given by
\begin{equation}
\begin{gathered}
x=\frac{x^1}{x^0}=\tan{\phi}, ~ ~
y=\frac{x^2}{x^0}=\tanh{r} \frac{\cos{(\theta-\phi)}}{\cos{\phi}}, \\
z=\frac{x^3}{x^0}=\tanh{r} \frac{\sin{(\theta-\phi)}}{\cos{\phi}}.
\end{gathered} \tag{2.10}
\end{equation}
\subsection{Geodesic curves}
The infinitesimal arc-length-square can be derived by the standard pull back method.
By $T^{-1}$-action of (2.6) on the differentials $(\mathrm{d}x^0;\mathrm{d}x^1;\mathrm{d}x^2;\mathrm{d}x^3)$, we obtain
that in this parametrization
the infinitesimal arc-length-square
at any point of $\SLR$ is the following:
\begin{equation}
   \begin{gathered}
      (\mathrm{d}s)^2=(\mathrm{d}r)^2+\cosh^2{r} \sinh^2{r}(\mathrm{d}\theta)^2+\big[(\mathrm{d}\phi)+\sinh^2{r}(\mathrm{d}\theta)\big]^2.
       \end{gathered} \tag{2.11}
     \end{equation}
Hence we get the symmetric metric tensor field $g_{ij}$ on $\SLR$ by components:
     \begin{equation}
       g_{ij}:=
       \begin{pmatrix}
         1&0&0 \\
         0&\sinh^2{r}(\sinh^2{r}+\cosh^2{r})& \sinh^2{r} \\
         0&\sinh^2{r}&1 \\
         \end{pmatrix}. \tag{2.12}
     \end{equation}

The geodesic curves of $\SLR$ are generally defined as having locally minimal arc length between any two of their (close enough) points.

By (2.9) the second order differential equation system of the $\SLR$ geodesic curve is the following:
\begin{equation}
\begin{gathered}
\ddot{r}=\sinh(2r)~\! \dot{\theta}~\! \dot{\phi}+\frac12 \big( \sinh(4r)-\sinh(2r) \big)\dot{\theta} ~\! \dot{\theta},\\
\ddot{\phi}=2\dot{r}\tanh{(r)}(2\sinh^2{(r)}~\! \dot{\theta}+ \dot{\phi}),\\ \ddot{\theta}=\frac{2\dot{r}}{\sinh{(2r)}}\big((3 \cosh{(2r)}-1)
\dot{\theta}+2\dot{\phi} \big). \tag{2.13}
\end{gathered}
\end{equation}

We can assume, by the homogeneity, that the starting point of a geodesic curve is the origin $(1,0,0,0)$.
Moreover, $r(0)=0,~ \phi(0)=0,~ \theta(0)=0,~ \dot{r}(0)=\cos(\alpha),~ \dot{\phi}(0)=\sin(\alpha)=-\dot{\theta}(0)$ are the initial values
in Table 1 for the solution of (2.12),
and so the unit velocity will be achieved.
\smallbreak
\centerline{\vbox{
\halign{\strut\vrule~\hfil $#$ \hfil~\vrule
&\quad \hfil $#$ \hfil~\vrule
&\quad \hfil $#$ \hfil\quad\vrule
&\quad \hfil $#$ \hfil\quad\vrule
&\quad \hfil $#$ \hfil\quad\vrule
\cr
\noalign{\vskip2pt}
\noalign{\hrule}
\multispan2{\strut\vrule\hfill{\bf  Table 1} \hfill\vrule}%
\cr
\noalign{\hrule}
\noalign{\vskip2pt}
{\rm Types} & {}  \cr
\noalign{\vskip2pt}
\noalign{\hrule}
\noalign{\vskip2pt}
\begin{gathered} 0 \le \alpha < \frac{\pi}{4} \\ (\bH^2-{\rm like~direction}) \end{gathered}
& \begin{gathered}  r(s,\alpha)={\mathrm{arsinh}} \Big( \frac{\cos{\alpha}}{\sqrt{\cos{2\alpha}}}\sinh(s\sqrt{\cos{2\alpha}}) \Big) \\
\theta(s,\alpha)=-{\mathrm{arctan}} \Big( \frac{\sin{\alpha}}{\sqrt{\cos{2\alpha}}}\tanh(s\sqrt{\cos{2\alpha}}) \Big) \\
\phi(s,\alpha)=2\sin{\alpha} s + \theta(s,\alpha) \end{gathered} \cr
\noalign{\vskip2pt}
\noalign{\hrule}
\noalign{\vskip2pt}
\begin{gathered} \alpha=\frac{\pi}{4} \\ ({\rm light~direction}) \end{gathered} &
\begin{gathered}  r(s,\alpha)={\mathrm{arsinh}} \Big( \frac{\sqrt{2}}{2} s \Big) \\
\theta(s,\alpha)=-{\mathrm{arctan}} \Big( \frac{\sqrt{2}}{2} s \Big) \\
\phi(s,\alpha)=\sqrt{2} s +\theta(s,\alpha) \end{gathered} \cr
\noalign{\vskip1pt}
\noalign{\hrule}
\noalign{\vskip2pt}
\begin{gathered} \frac{\pi}{4}  < \alpha \le \frac{\pi}{2} \\ ({\rm fibre-like~direction}) \end{gathered} &
\begin{gathered}  r(s,\alpha)={\mathrm{arsinh}} \Big( \frac{\cos{\alpha}}{\sqrt{-\cos{2\alpha}}}\sin(s\sqrt{-\cos{2\alpha}}) \Big) \\
\theta(s,\alpha)=-{\mathrm{arctan}} \Big( \frac{\sin{\alpha}}{\sqrt{-\cos{2\alpha}}}\tan(s\sqrt{-\cos{2\alpha}}) \Big) \\
\phi(s,\alpha)=2\sin{\alpha} s + \theta(s,\alpha) \end{gathered}  \cr
\noalign{\vskip2pt}
\noalign{\hrule}
\noalign{\hrule}}}}
\smallbreak
The equation of the geodesic curve in the hyperboloid model can be given by the usual geographical sphere coordinates
$(\lambda, \alpha)$, as longitude and altitude, respectively, from the general starting position of (2.10), (2.11),
$(-\pi < \lambda \le \pi, ~ -\frac{\pi}{2}\le \alpha \le \frac{\pi}{2})$,
and the arc-length parameter $0 \le s \in \bR$. The Euclidean coordinates $X(s,\lambda,\alpha)$,
$Y(s,\lambda,\alpha)$, $Z(s,\lambda,\alpha)$ of the geodesic curves can be determined by substituting the results of Table 1 into the
equations (2.11) as follows (see e.g. \cite{MSz})
\begin{equation}
\begin{gathered}
X(s,\lambda,\alpha)=\tan{(\phi(s,\alpha))}, \\
Y(s,\lambda,\alpha)=\frac{\tanh{(r(s,\alpha))}}{\cos{(\phi(s,\alpha)}} \cos \big[ \theta(s,\alpha)-\phi(s,\alpha)+\lambda \big],\\
Z(s,\lambda,\alpha)=\frac{\tanh{(r(s,\alpha))}}{\cos{(\phi(s,\alpha)}} \sin \big[ \theta(s,\alpha)-\phi(s,\alpha)+\lambda \big].
\end{gathered} \tag{2.14}
\end{equation}
\begin{Definition}
The {\rm geodesic distance} $d(P_1,P_2)$ between the points $P_1$ and $P_2$ is defined by the arc length of the geodesic curve
from $P_1$ to $P_2$.
\end{Definition}
\section{Generalized $\SLR$ cylinders}
Let $\cC_c^b(r)$ denote the circle of $\SLR$ radius $r \in \mathbf{R}^+$ 
centered at the origin, lying in the base plane of the model and  
let $\cC^b(\infty)$ denote the base circle 
of the hyperboloid model, centred at the origin and lying in the base plane.
\begin{Definition}
\begin{enumerate}
\item Let $\cC^i_c(r)$ be an infinite solid that is ruled by ``side fibre lines" passing 
through the points of a circle $\cC^b_c(r)$ of $\SLR$ radius $r\in \mathbf{R}$ $(r \ge 0)$ lying in 
the base plane and centred at the origin.
The images $\cC^i_{c,\bt(r)}$ of solid $\cC^i_c(r)$ by $\SLR$ isometries $\bt$ 
are called {\rm infinite circular cylinders}. (We also include the fiber lines in this class if $r=0$).
\item  Let $\cC^i_h$ be a surface that is generated by ``side fibre lines" passing 
through the points of arc of a Euclidean circle $\cC^b_h$ lying in 
the base plane where $\mathrm{card}\{\cC^i_h \cap \cC^b(\infty)\}=2$ 
The images $\cC^i_{h,\bt}$ of solid $\cC^i_h$ by $\SLR$ isometries $\bt$ 
are called {\rm infinite hypercyclical cylinders}.
\item Let $\cC^i_p$ be a surface that is generated by ``side fibre lines" passing 
through the points an Euclidean circle $\cC^b_p$ lying in the closed region bounded by $\cC^b(\infty)$
in the the base plane where $\mathrm{card}\{\cC^i_p \cap \cC^b(\infty)\}=1$.  
The images $\cC^i_{p,\bt}$ of solid $\cC^i_p$ by $\SLR$ isometries $\bt$ 
are called {\rm infinite horocyclical cylinders}.
\item Surfaces that satisfy the previous 
three conditions are uniformly called {\rm generalized circular cylinders} in $\SLR$ geometry.
\end{enumerate}
\end{Definition}
In \cite{Sz24} we proved the following lemmas for $\cC^i_c(r)$ cylinders (for classical cylinder types,
see Definition 3.1 (1)).
\begin{lemma}[\cite{Sz24}]
If the radius of the circle $\cC^b_c(r)$ centred at the origin lying in the base plane 
is $r$ then the corresponding infinite fibre-like cylinder $\cC_c^i(r)$
is a one-sheeted hyperboloid in Euclidean sense with equation
\begin{equation}
\frac{y^2}{\tanh^2{r}}+\frac{z^2}{\tanh^2{r}}-x^2=1, \quad {\text{where}} 
\quad x=\frac{x^1}{x^0},~y=\frac{x^2}{x^0},~z=\frac{x^3}{x^0}. \tag{3.1}
\end{equation}
\end{lemma}
Let us denote the image of the cylinder $\cC^i_c(r)$ at the $\SLR$ translation $\bT$ (see (2.5))
by $\cC^i_{c,\bT}(r)$ where the translation is given by parameters $x^0=1,x^1,x^2,x^3$ and denote
its common part with the base plane by $\cC^b_{c,\bT}(r)$.
\begin{lemma}[\cite{Sz24}]
$\cC^b_{c,\bT}(r)$ is a circle of radius
\begin{equation}
R=\frac{\tanh{r}(1+(x^1)^2-(x^2)^2-(x^3)^2)}{1+(x^1)^2-\tanh^2{r}\cdot((x^2)^2+(x^3)^2)}, \tag{3.2}
\end{equation}
in Euclidean sense and its centre can be derived by the point $K=(1;0;\frac{x^2-x^1x^3}{1+x^1x^1};$ 
$\frac{x^3+x^1x^2}{1+x^1x^1})$ and by its radius $r$. 
\end{lemma}
The previous lemma can be formulated and proved in a similar way for 
other types of generalized cylinders (see Definition 3.1 (2),(3)).
Let us denote the image of a horocyclical or hypercyclical cylinder $\cC^i_j(r)$ 
$(j \in \{h,p\})$ at the $\SLR$ translation $\bT$ (see (2.5))
by $\cC^i_{j,\bT}(r)$ where the translation is given by $(x^0=1;x^1;x^2;x^3)$ and denote
its common part with the base plane by $\cC^b_{j,\bT}(r)$.
\begin{lemma}
$\cC^b_{j,\bT}(r)$ $(j \in \{h,p\})$ is a circular arc 
if the cylinder is hypercyclical or a circle if the cylinder is horocyclical 
that touches $\cC^b(\infty)$, the base circle of the model. \quad \quad $\square$
\end{lemma}
\begin{Remark}
\begin{enumerate}
\item The translations correspond a fibrum to a fibrum, so it is sufficient to consider 
the translations in the base plane, that means can be assumed $x^1=0$.
\item The perimeter $p(\cC^b_c(r))$ of a circle $\cC^b_c(r)$ and the area 
${\mathrm{Area}}(\cC^b_c(r))$ 
with radius $r$ ($\SLR$ geodesic distance) 
lying in the base plane centred at the origin of the projective model of the $\SLR$
geometry were determined in \cite{Sz24}. 
\item Moreover,in \cite{Sz24} we determined the surface area ${\mathrm{SA}}(\cC(r,\Psi))$ and 
the volume ${\mathrm{Vol}}((\cC(r,\Psi)))$ of a bounded fibre-like circular cylinder with radius $r$.
\end{enumerate}
\end{Remark}
\section{Relationship between cylinders and geodesic curves in $\SLR$ geometry}
In Section 2, we provided the equations for the types of geodesic curves 
in Table 1. The images of these types, as determined by the $\SLR$ isometries, 
constitute the equivalence classes of the geodesics.
As can be seen from the table, the equivalence classes are determined 
by the value of the parameter $\alpha$ ($0 \le \alpha \le \frac{\pi}{2}$).
\begin{lemma}
The geodesic curves lie on generalized cylindrical surfaces (see Definition 3.1):
\begin{enumerate}
\item If $0 \le \alpha \le \frac{\pi}{4}$ then on {\rm infinite hypercyclical cylinders},
\item If $\alpha=\frac{\pi}{4}$ then on {\rm infinite horocyclical cylinders}.
\item If $\frac{\pi}{4} \le \alpha \le \frac{\pi}{2}$ on {\rm infinite circular cylinders}.
\end{enumerate}
where the Euclidean radii of the intersections - circles or circular arcs - of the generalized cylinders with a base plane are
$R=\frac{1}{2 \cdot \tan{\alpha}}$ (if $\alpha=\frac{\pi}{2}$ then $R:=0$).
\end{lemma}
{\bf Proof:}

The fibre line through $X=(x^0;x^1;x^2;x^3)$ intersects the base plane in the foot point $X'$ (see (2.3), (2.4), (2.8))
\begin{equation}
\begin{gathered}
X'=\Big(1;0; \frac{x^0x^2-x^1x^3}{x^0 x^0+x^1x^1};\frac{x^0x^3+x^1x^2}{x^0 x^0+x^1x^1}\Big).
\end{gathered} \tag{4.1}
\end{equation}
If we apply formula (4.1) to the points of the geodesic curve (see Table 1), we obtain a one-parameter planar curve 
in the base plane of the model; determining its curvature can be seen that the radius of curvature depends solely on the parameter $\alpha$  
and its value is $\frac{1}{2 \cdot \tan{\alpha}}$
implying that the projection of the geodesic curve is a circle or a circular arc within the model's base circle in the Euclidean sense.

Based on Lemmas 3.2 and 3.3, this means that for every possible parameter $\alpha$, the geodesic curves lie on generalized cylindrical surfaces. \quad \quad $\square$

\subsection{Pappus-like theorem for circles in inversive plane}
The classical M\"obius plane is the Euclidean plane supplemented by a single point at infinity. It is also called the inversive plane because it is closed under 
inversion with respect to any circle. The following theorem can be proven using the usual inversion and elementary geometric theorems and relations. 
\begin{Theorem}
Consider six arbitrary points on the inversive plane 
$A_1, A_2, A_3,$ $B_1, B_2, B_3$. We consider two circles (can also be a straight line, i.e. a circle with infinite radius) 
$K_{A_1A_2A_3}$ passing through the points $A_1,A_2,A_3$ and $K_{B_1B_2B_3}$ containing $B_1, B_2, B_3$. Let a point $\overline{A}_1$ be on the circle $K_{A_1A_2A_3}$. 
Let $\overline{B}_2$ be the point where the circle $K_{A_1B_2\overline{A}_1}$ intersects the circle $K_{B_1B_2B_3}$. Moreover, draw the circle $K_{B_2A_3\overline{B}_2}$ and let 
$\overline{A}_3=K_{B_2A_3\overline{B}_2} \cap K_{A_1A_2A_3}$ where $\overline{A}_3 \ne \overline{B_2}$. Let $\overline{B}_1$ be the from $B_1$ different 
point where the circle $K_{A_3B_1\overline{A}_3}$
intersects the circle $K_{B_1B_2B_3}$. Then, we consider the circle $K_{B_1A_2\overline{B}_1}$ and let $\overline{A}_2$ be the other point where it intersects the circle $K_{A_1A_2A_3}$. 
Let $\overline{B}_3$ be the intersection point of the circles $K_{A_2B_3\overline{A}_2}$ and $K_{B_1B_2B_3}$ which is different from $B_3$. 
\begin{figure}[ht]
\centering
\includegraphics[width=13cm]{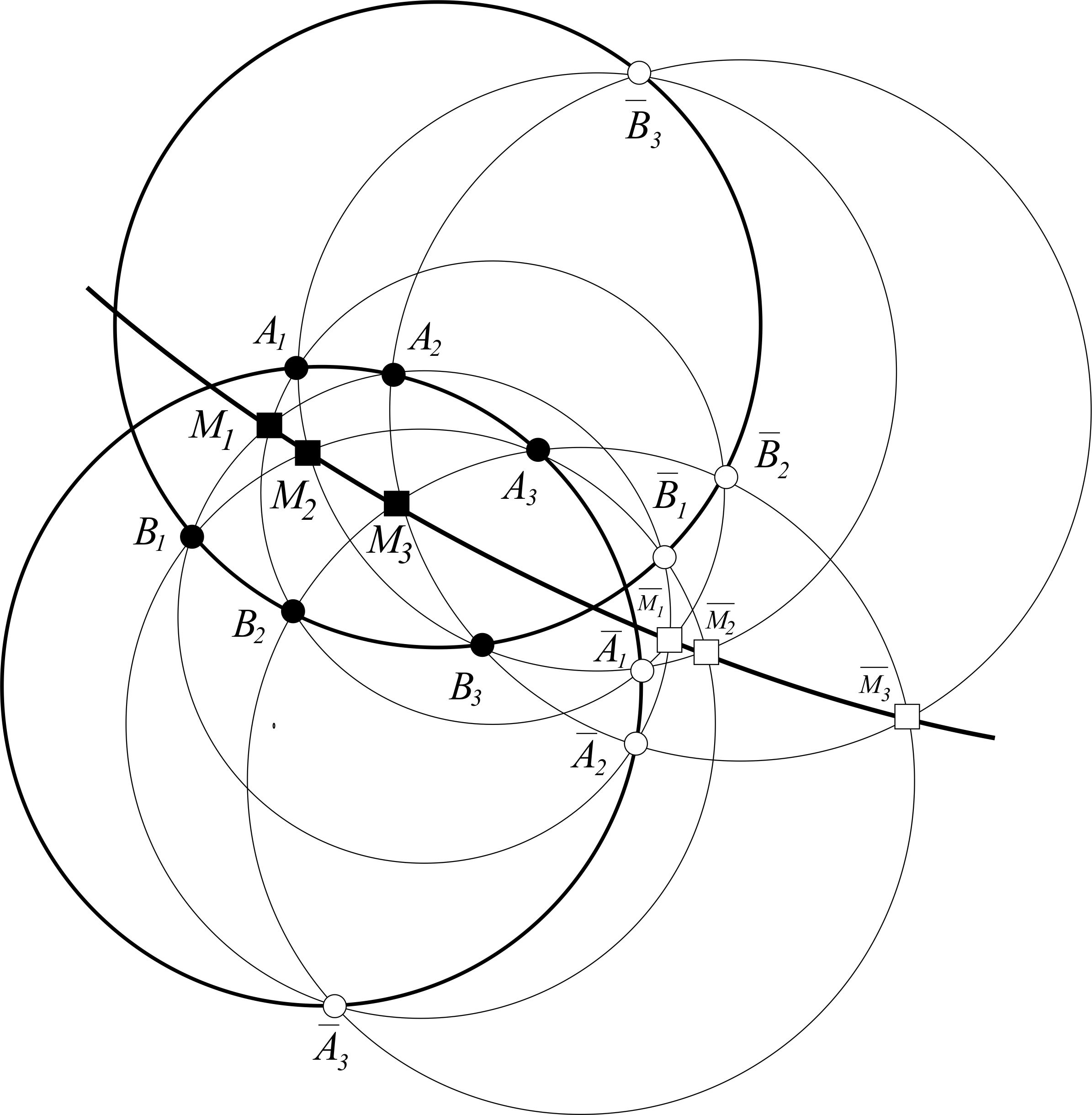}
\caption{Circle-like Pappus's hexagon theorem in the inversive plane}
\label{}
\end{figure}
\begin{enumerate}
\item Thus, $\overline{A}_1 \in K_{B_3A_1\overline{B}_3}$. 
\item Let $M_1$ and $\overline{M}_1$ be the points where the circles $K_{A_1B_2 \overline{A}_1 \overline{B}_2}$ and $K_{A_2B_1 \overline{A}_2 \overline{B}_1}$ meet, 
and $M_2$ and $\overline{M}_2$ where 
$K_{A_1B_3\overline{A}_1,\overline{B}_3}$ and $K_{A_3B_1\overline{A}_3,\overline{B}_1}$ meet, moreover, $M_3$ and $\overline{M}_3$ 
where $K_{A_3B_2\overline{A}_3,\overline{B}_2}$ and $K_{A_2B_3\overline{A}_2,\overline{B}_3}$ meet. 
Then $M_1, M_2, M_3$, $\overline{M}_1$, $\overline{M}_2$ and $\overline{M}_3$ are concyclic. 
\end{enumerate} \qquad $\square$
\end{Theorem}
\begin{Remark}
If the points $\overline{A}_i$ and $\overline{B}_i$ $i \in \{1,2,3 \}$ coincide, then using an
inversion we directly obtain the classical Pappus hexagon configuration.
\end{Remark}
\subsection{Pappus-like theorem for generalized $\SLR$ cylinders}
Let us formulate the next  theorem similar to Theorem 4.1 for $\SLR$ generalized cylinders.
\begin{figure}[ht]
\centering
\includegraphics[width=12cm]{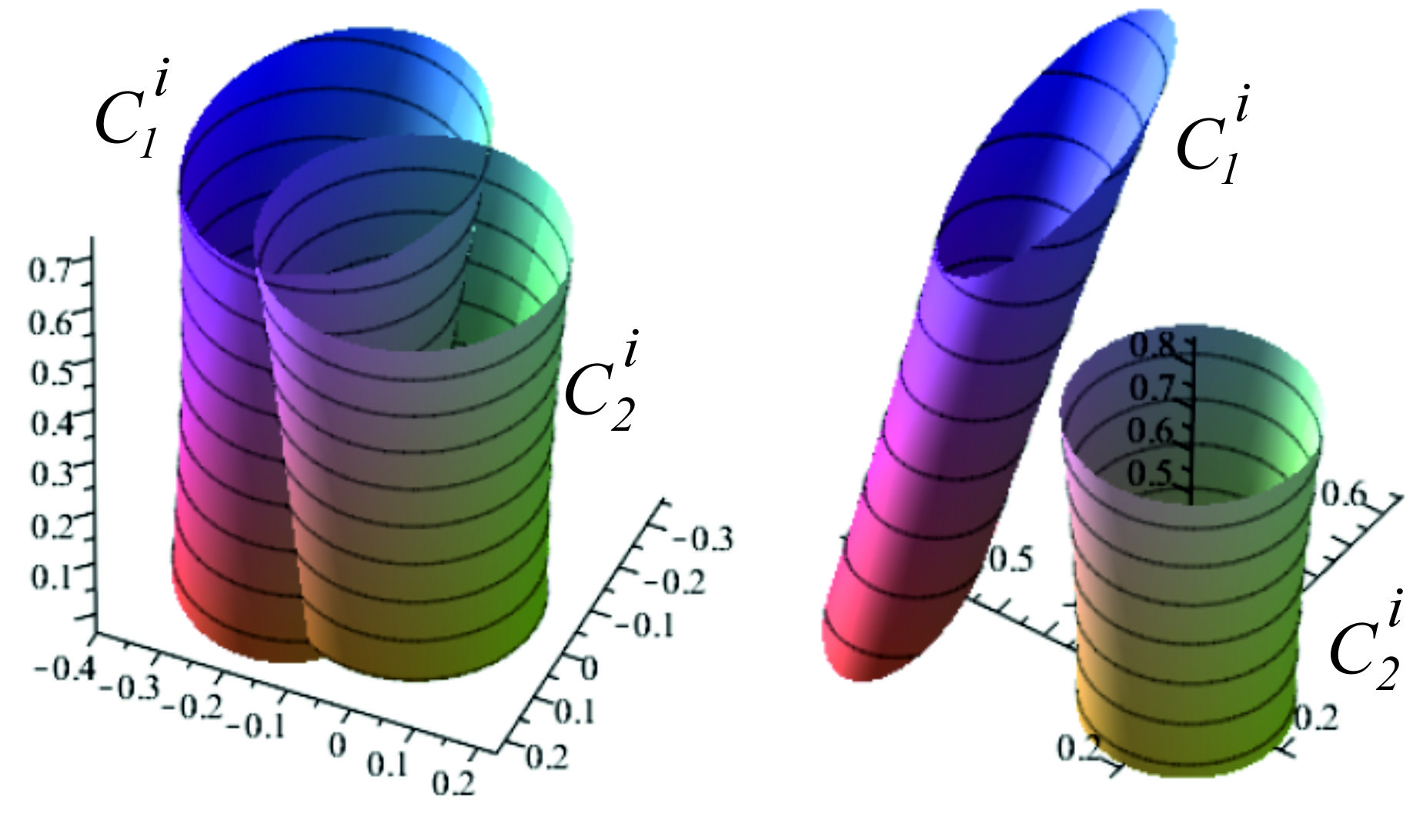}
\caption{Two intersecting congruent circular cylinders (left). 
\ \ Two congruent circular cylinders (right).}
\label{}
\end{figure}
\begin{Theorem}
Consider six arbitrary distinct points on the base plane of the considered model of the $\SLR$ geometry  
$A_1, A_2, A_3$ and $B_1, B_2, B_3$ that circles defined by which lie on generalized cylinders. 

We assume that the intersection points arising in the base plane - which appears in the theorem - as well as the corresponding fibres lie in the interior or on the boundary  
of the considered model of $\SLR$ space.

We consider two {\it generalized cylinders} (can also be a Euclidean plane, i.e. a hypercyclical cylinder with infinite radius) 
$\cC^{i}_{A_1A_2A_3}$ passing through the points $A_1,A_2,A_3$ and $\cC^i_{B_1B_2B_3}$ containing $B_1, B_2, B_3$. 
Let a fibre line $\overline{a}^f_1$ be on the cylinder 
$\cC^{i}_{A_1A_2A_3}$ and let its common point be with the base plane $\overline{A}_1$ 
Let $\overline{b}^f_2$ be the fibre line where the cylinder $\cC^i_{j,A_1B_2\overline{A}_1}$ intersects the cylinder $\cC^i_{B_1B_2B_3}$ and we denote 
its intersection point with the base plane by $\overline{B}_2$. Moreover, consider the cylinder $\cC_{B_2A_3\overline{B}_2}$ and let 
$\overline{a}^f_3=\cC^i_{B_2A_3\overline{B}_2} \cap \cC^i_{A_1A_2A_3}$ where $\overline{A}_3 \ne \overline{B}_2$. Let $\overline{b}^f_1$ be a fibre line 
where the cylinder $\cC^i_{A_3B_1\overline{A}_3}$
intersects the cylinder $\cC^i_{B_1B_2B_3}$ moreover its common point $\overline{B}_1$ with the base plane is different from $B_1$. 
Then, we consider the cylinder $\cC^i_{B_1A_2\overline{B}_1}$ and let $\overline{a}^f_2$ be the other fibre line where it intersects the cylinder $\cC^i_{A_1A_2A_3}$ and 
let $\overline{A}_2$ the intersection pont of $\overline{a}^f_2$ with the base plane of the model. 
Let $\overline{b}^f_3$ be the intersection fibre line of the cylinders $\cC^i_{A_2B_3\overline{A}_2}$ and $\cC^i_{B_1B_2B_3}$ which is different from the fibre line through $B_3$, 
moreover, we denote its intersection point with the base plane by $\overline{B}_3$.
\begin{figure}[ht]
\centering
\includegraphics[width=12cm]{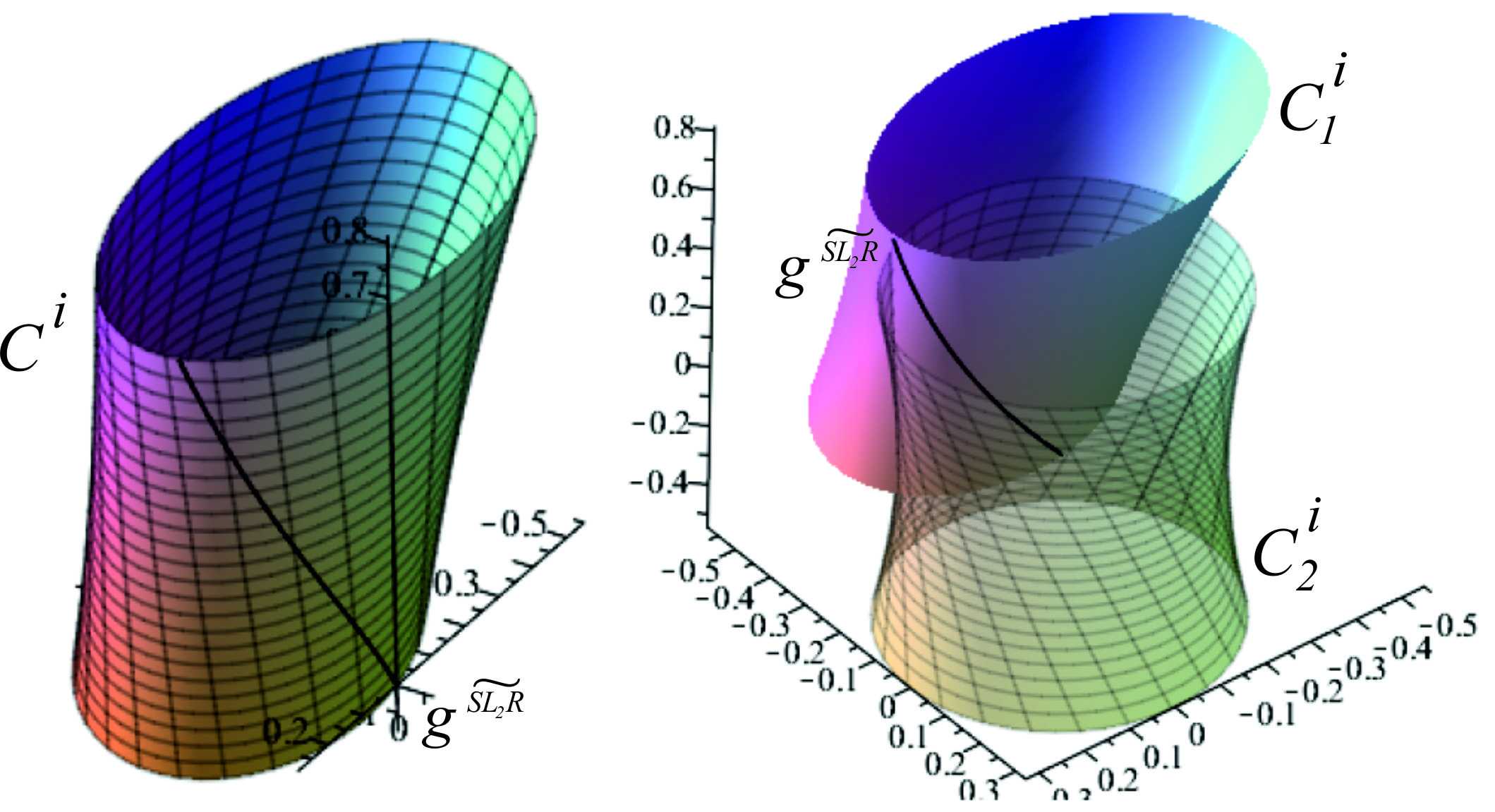}
\caption{Geodesic curve with the circular cylinder containing it (left). \ \ Two congruent circular cylinders with a geodesic curve (right).}
\label{}
\end{figure}
\begin{enumerate}
\item Thus, $\overline{a}^f_1 \in \cC^i_{B_3A_1\overline{B}_3}$. 
\item Let $m^f_1$ and $\overline{m}^f_1$ be the fibre lines where the cylinders $\cC^i_{A_1B_2 \overline{A}_1 \overline{B}_2}$ and $\cC^i_{A_2B_1 \overline{A}_2 \overline{B}_1}$ meet, 
$m^f_2$ and $\overline{m}^f_2$ where 
$\cC^i_{A_1B_3\overline{A}_1,\overline{B}_3}$ and $\cC^i_{A_3B_1\overline{A}_3,\overline{B}_1}$ meet, moreover, $m_3^f$ and $\overline{m}_3^f$ 
where $\cC^i_{A_3B_2\overline{A}_3,\overline{B}_2}$ and $\cC^i_{A_2B_3\overline{A}_2,\overline{B}_3}$ meet. 
Then the fibre lines $m_1^f, m_2^f, m_3^f$, $\overline{m}_1^f, \overline{m}_2^f$ and $\overline{m}_3^f$ lie in a generalized cylinder. 
\end{enumerate} \qquad $\square$
\end{Theorem}
Based on the Definition 3.1, Lemmas 3.1, 3.2, 3.3, 4.1 the types of generalized cylinders determine equivalence classes of geodesic curves up to $\SLR$ congruences, 
i.e., the Euclidean radii of the 
$\SLR$ infinite generalized cylinders characterizes the equivalence classes of geodesic curves. 
Based on these, using the notations of Theorem 4.2 we can state the following theorems:
\begin{Theorem}
There exists uniquely equivalence class of $\SLR$ geodesics that intersects the fibres $m_1^f, m_2^f, m_3^f, \overline{m}_1^f, \overline{m}_2^f$ 
and $\overline{m}_3^f$ \qquad \qquad $\square$
\end{Theorem}
It would not be easy to visualize the many cylinders that occur in the theorem, 
so we only show the intersection of a geodesic with a fibre line in Fig.~4.
\begin{figure}[ht]
\centering
\includegraphics[width=12cm]{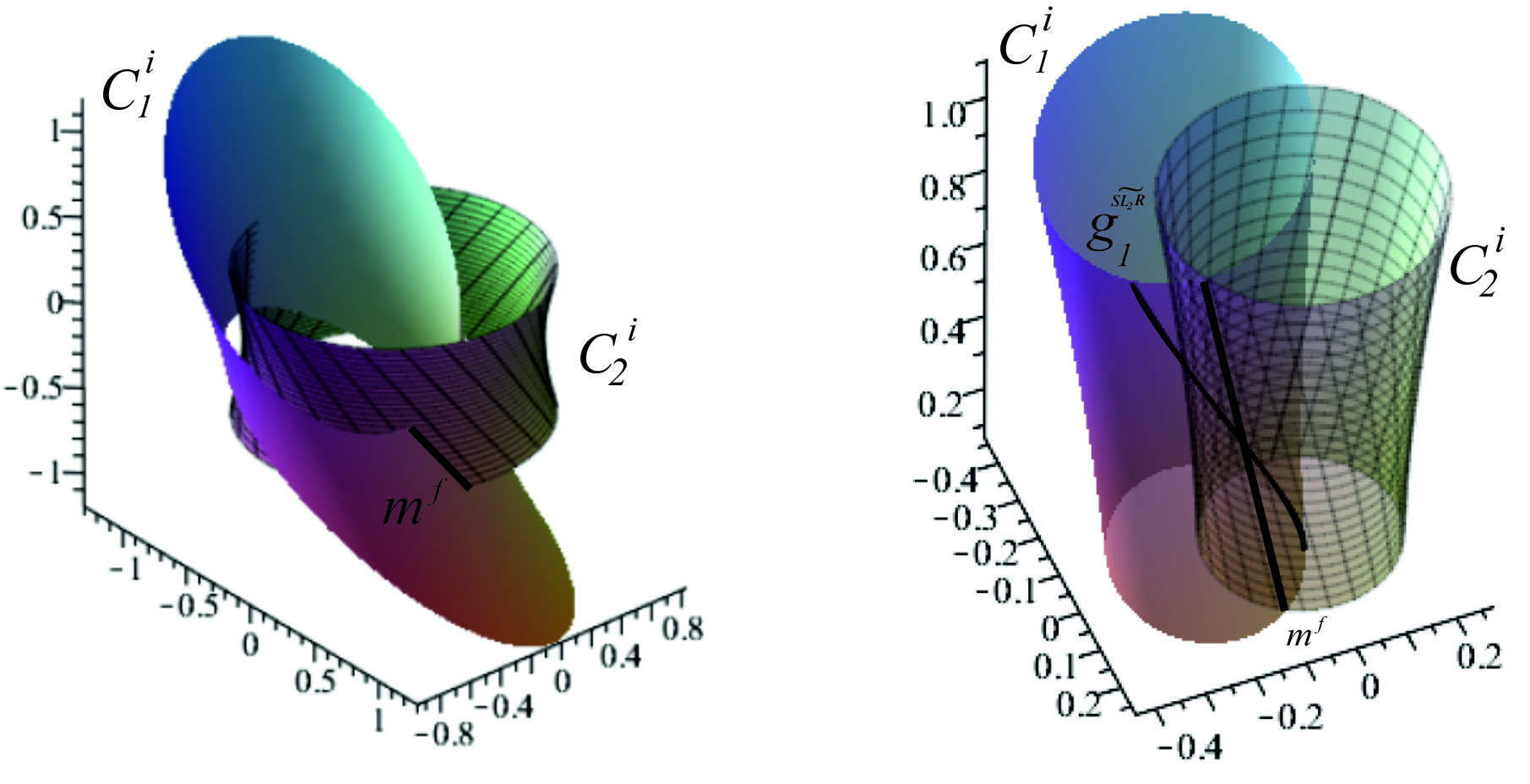}
\caption{Two circular cylinders with a common fibre line (left). \ \ Two circular cylinders and the intersection of a geodesic curve with a fibre line (right).}
\label{}
\end{figure}
\begin{Theorem}
Let $g^{\SLR}_{m_1^fm_2^fm_3^f}\subset 
\cC^i_{m_1^fm_2^fm_3^f}$ be an given geodesic curve (see Theorem 4.3). 
There always exist representation elements of the equivalence classes of $\SLR$ geodesics related to the Theorem 4.3
so that their intersection points 
\begin{enumerate}
\item $M_1 = g^{\SLR}_{A_1B_2 \overline{A}_1 \overline{B}_2} \cap g^{\SLR}_{A_2B_1 \overline{A}_2 \overline{B}_1} \in m_1^f$,
\item $M_2 = g^{\SLR}_{A_1B_3\overline{A}_1,\overline{B}_3} \cap g^{\SLR}_{A_3B_1\overline{A}_3,\overline{B}_1} \in m_2^f$,
\item $M_3 = g^{\SLR}_{A_3B_2\overline{A}_3,\overline{B}_2} \cap g^{\SLR}_{A_2B_3\overline{A}_2,\overline{B}_3} \in m_3^f$
\end{enumerate}
exist and lie on $g^{\SLR}_{m_1^fm_2^fm_3^f}\subset \cC^i_{m_1^fm_2^fm_3^f}$. \qquad \qquad $\square$
\end{Theorem}
As it is known, the Pappus's hexagon theorem plays a key role in the axiomatic structure of projective geometry. We 
proved that an analogous theorem is also true in $\NIL$ (see \cite{Sz26}) space and we have seen that in $\SLR$ geometry, too. Therefore, it is interesting to 
examine the role of this theorem in the structure of in other non-constant curvature Thurston geometries.
\medbreak
In this paper we mentioned only some natural problems related to 
Thurston spaces, but we hope that from these
the reader can appreciate that our projective method is suitable to study and solve 
similar problems that represent a huge class of open mathematical problems 
(see e.g. \cite{Cs23}, \cite{CsSz23}, \cite{MSz12}, \cite{M-Sz}, \cite{stachel}, \cite{MSzV13},
\cite{Sz07}, \cite{Sz10-2}, \cite{Sz10-3}, \cite{Sz12-1}, \cite{Sz13-1}, 
\cite{Sz14}, \cite{Sz13-2}, 
\cite{Sz202}, \cite{Sz26-0}, \cite{Sz26-1}).
Detailed studies are the objective of ongoing research.

\vspace{0.5cm}

{\bf Funding}: The authors declare that no funds, grants, or other support were
received during the preparation of this manuscript.

\vspace{2mm}

{\bf Data Availability Statement}: The authors declare that this research does not
use data from any external sources or counterparts.

\vspace{2mm}

{\bf Declarations}: I declare that this paper presented in it is my own and has 
been generated by me as the result of my own original research.

\vspace{2mm}

{\bf Competing Interest}: The authors have no relevant financial or non-financial
interests to disclose.


\begin{thebibliography}{MPSz98}
%
\bibitem{Cs23}
{Csima,~G.:}
Isoptic surfaces of segments in $\SXR$ and
$\HXR$ geometries.
\emph{J. Geom.} (2024), doi: 10.1007/s00022-023-00699-x.
%
\bibitem{CsSz16}
Csima,~G.~--~Szirmai,~J.:
Interior angle sum of translation and geodesic triangles in $\SLR$ space, 
{\it Filomat}, {\bf 32/14} (2018), 5023--5036. 
%
\bibitem{CsSz23}
Csima,~G.~--~Szirmai,~J.:
Translation-like isoptic surfaces and angle sums of translation triangles in $\NIL$ geometry, 
{\it Results Math.}, (2023) {\bf 78:194}, doi: 10.1007/s00025-023-01961-z. 
%
\bibitem{M97}
{Moln{\'a}r,~E.}
The projective interpretation of the eight 3-di\-men\-sional homogeneous geometries.
\emph{Beitr. Algebra Geom.,}
{\bf38} (1997), No.~2, 261--288.
  %
\bibitem{MSz}
{Moln{\'a}r,~E.~--~Szirmai,~J.}
Symmetries in the 8 homogeneous 3-geometries.
\emph{Symmetry Cult. Sci.,}
{\bf 21/1-3} (2010), 87-117.
%
\bibitem{MSz12}
{Moln{\'a}r,~E.~--~Szirmai,~J.}
Classification of $\SOL$ lattices.
\emph{Geom. Dedicata,}
{\bf 161/1} (2012), 251-275, doi: 10.1007/s10711-012-9705-5.
%
\bibitem{MSz14}
{Moln{\'a}r,~E.~--~Szirmai,~J.}
Volumes and geodesic ball packings to the regular prism tilings in $\SLR$ space.
\emph{Publ. Math. Debrecen,}
{\bf 84/1-2} (2014), 189-203, doi: 10.5486/PMD.2014.5832.
%
\bibitem{M-Sz}
{Moln{\'a}r,~E.~--~Szirmai,~J.}
On homogeneous 3-geometries, balls and their optimal arrangements, especially 
in $\NIL$ and $\SOL$ spaces.
\emph{G-Slovak Journal for Geometry and Graphics} {\bf 19}  (2022), 37, 5-32.
%
\bibitem{stachel}
Moln{\'a}r,~E.~---~Szirmai,~J.
Packings with geodesic and translation balls and their visualizations in $\SLR$ space.
{\it Journal for Geometry and Graphics}  {\bf 26}(1) (2022), 51--64.
%
\bibitem{MSzV}
 {Moln{\'a}r,~E.~--~Szirmai,~J.~--~Vesnin,~A.}
 Projective metric realizations of cone-manifolds with singularities along 2-bridge knots and links.
 \emph{J. Geometry,} {\bf 95} (2009), 91-133.
%
\bibitem{MSzV13}
 {Moln{\'a}r,~E.~--~Szirmai,~J.~--~Vesnin,~A.}
 Packings by translation balls in $\SLR$. \emph{J. Geometry,}
 {\bf 105 (2)} (2014), 287-306, doi: 10.1007/s00022-013-0207-x.
%
\bibitem{R}
{Ratcliffe,~J.~G.}
\emph{Foundations of hyperbolic manifolds, (2nd ed.).}
{Graduate Texts in Mathematics 149. New York, NY: Springer.,} (2006).
 %
\bibitem{S}
{Scott,~P.}
The geometries of 3-manifolds.
\emph{Bull. London Math. Soc.}, {\bf15} (1983), 401--487.
  %
\bibitem{Sz07}
{Szirmai,~J.}
The densest geodesic ball packing by a type of $\NIL$ lattices.
\emph{Beitr. Algebra Geom.,}
 {\bf48(2)} (2007), 383--398.
 %
\bibitem{Sz10-2}
{Szirmai,~J.}
Geodesic ball packing in $\SXR$ space for generalized Coxeter space groups.
\emph{Beitr. Algebra Geom.,}
{\bf 52(2)} (2011), 413--430.
 %
\bibitem{Sz10-3}
{Szirmai,~J.}
Geodesic ball packing in $\HXR$ space for generalized Coxeter space groups.
\emph{Math. Commun.,}
{\bf 17/1} (2012), 151-170.
 %
\bibitem{Sz12-1}
{Szirmai,~J.}
Lattice-like translation ball packings in $\NIL$ space.
\emph{Publ. Math. Debrecen,}
{\bf 80/3-4} (2012), 427--440 doi: 10.5486/PMD.2012.5117.
  %
\bibitem{Sz13-1}
{Szirmai,~J.}
Regular prism tilings in $\SLR$ space.
\emph{Aequat. Math.}, (2013), doi 10.1007/s00010-013-0221-y.
%
\bibitem{Sz14}
Szirmai,~J
Non-periodic geodesic ball packings to infinite regular prism tilings in $\SLR$ space.  
\emph{Rocky Mountain J.  Math.}  {\bf 46/3} (2016), 1055--1070.
%
\bibitem{Sz13-2}
{Szirmai,~J.}
A candidate to the densest packing with equal balls in the Thurston geometries.
\emph{Beitr. Algebra Geom.,} (2013), doi: 10.1007/s13366-013-0158-2.
%
\bibitem{Sz202} Szirmai,~J.,
Interior angle sums of geodesic triangles in $\SXR$ and $\HXR$ geometries, 
{\it Bull. Academ. De Stiinte A Rep. Mol.}, {\bf{93}} Num 2 (2020), 44--61.
%
\bibitem{Sz22-2}
Szirmai,~J.,
Apollonius surfaces, circumscribed spheres of tetrahedra, Menelaus' and Ceva's theorems in $\SXR$ and $\HXR$ geometries. 
{\it Q. J. Math.,} {\bf 73}, (2022), 477-494, doi: 10.1093/qmath/haab038.
%
\bibitem{Sz23-1}
Szirmai,~J.,
On Menelaus' and Ceva's theorem in $\NIL$ geometry. 
{\it Acta Univ. Sapientiae Math.,} {\bf 15/ 1}, (2023) 123--141, doi: 10.2478/ausm-2023-0008.
%
\bibitem{Sz23-2}
{Szirmai,~J.,}
Classical Notions and Problems in Thurston Geometries,
\emph{International Electronic Journal of Geometry},
{\bf 16} No.2 (2023), 608--643, doi: 10.36890/IEJG.1221802.
%
\bibitem{Sz24}
{Szirmai,~J.},
Fibre-like cylinders, their packings and coverings in $\SLR$ space,
\emph{Results Math.}, {\bf{79}} article number 123 (2024), doi: 10.1007/s00025-024-02152-0.
%
\bibitem{Sz26-0}
{Szirmai,~J.},
Menelaus' and Ceva's theorems for translation triangles in Thurston geometries, 
\emph{Results Math.}, (2026), {\bf 81:160}, DOI: 10.1007/s00025-026-02706-4.
%
\bibitem{Sz26-1}
{Szirmai,~J.},
Desargues's and Pappus's hexagon theorems on translation-type surfaces in Thurston geometries, 
\emph{Aequat. Math.}, (2026), {\bf 100:82}, DOI: 10.1007/s00010-026-01325-2.
%
\bibitem{Sz26}
{Szirmai,~J.:}
Cylinder-like Pappus's hexagon theorem in $\NIL$ geometry,
\emph{Submitted manuscript}, (2026), arXiv: 2606.27765.
%
\bibitem{T}
{Thurston,~W.~P.} (and {\sc Levy,~S.} editor)
\emph{Three-Dimensional Geometry and Topology.}
{Princeton University Press,} Princeton, New Jersey, Vol.{\bf 1} (1997).

 \end{thebibliography}
 \end{document}